\documentclass[12pt]{amsart}

\usepackage{amsfonts,amssymb,amsmath,amsgen} 
\usepackage{appendix}
\newtheorem{theorem}{Theorem}[section]

\newtheorem{lemma}[theorem]{Lemma}
\newtheorem{proposition}[theorem]{Proposition}
 \newtheorem{remark}[theorem]{Remark}
 \newtheorem{definition}[theorem]{Definition}

\newcommand{\R}{  \mathbb{R}   }
\newcommand{\eps}{\epsilon}
\newcommand{\e}{\epsilon}

\begin{document}
 \bibliographystyle{plain}
\title{Unconditional well-posedness for wave maps}

   \author{Nader Masmoudi}
\address{ Courant Institute of Mathematical Sciences\\
251 Mercer Street, New York NY 10012, U.S.A.\\
masmoudi@cims.nyu.edu}

 \author{Fabrice Planchon}
 \address{Laboratoire J. A. Dieudonn\'e, UMR 6621\\
 Universit\'e de Nice Sophia-Antipolis\\
 Parc Valrose\\
 06108 Nice Cedex 02, FRANCE\\
et Institut universitaire de France\\
fabrice.planchon@unice.fr}

 \date{} \maketitle
 \begin{abstract}
We prove uniqueness of solutions to the wave map equation
  in the natural class, namely  $ (u, \partial_t u ) \in  C([0,T); \dot{H}^{d/2})\times
C^1([0,T); \dot{H}^{d/2-1})$   in dimensions $d\geq 4$. This is
achieved through estimating the difference of two solutions at a lower
regularity level. In order to reduce to the Coulomb gauge, one has to
localize the gauge change in suitable cones as well as estimate the
difference between the frames and connections associated to each
solutions and take advantage of the assumption that the target
manifold has bounded curvature.
 \end{abstract}
 \par \noindent

\section{Introduction}
Let $(N,g)$ be  a  complete riemannian manifold of dimension
$k$ without boundary. We denote  $(x^\alpha), \, 0\leq \alpha \leq d$ the
canonical coordinate system of $ \mathbb{R} \times \mathbb{R}^d $ where
$t = x^0$ denotes the time variable. Moreover, we denote $\partial_\alpha
= \partial/\partial x^\alpha$ and use the Minkowski metric
on $ \mathbb{R} \times \mathbb{R}^d $ to raise and lower indices.
In particular, $ \partial^0 = - \partial_0 $ and
$ \partial^\alpha =  \partial_\alpha  $ for $1\leq \alpha \leq d$.
The wave map  equation, from  $ \mathbb{R} \times \mathbb{R}^d $
into $N$,   reads

\begin{equation}
  \label{wm}
\left \{ \begin{array}{rcl} D_\alpha \displaystyle{\partial}^\alpha u
    & = & 0,\\
u(x,0) & = & u_0(x),\\
\partial_t u(x,0) & = & u_1(x) \quad \quad
  x\in \mathbb{R}^d,\quad t\geq 0\,,
\end{array} \right.
\end{equation}
where $D_\alpha$ is the pull-back of the
 covariant derivative on the target Riemannian
manifold $N$.

\subsection{Statement of the  Result}
Low regularity strong (e.g. unique) solutions to semilinear wave
equations like (\ref{wm}) are usually constructed via
fixed point methods. Hence, while one is ultimately seeking solutions
which are continuous evolutions of the data, that is $(u,\partial_t u)
\in C([0,T); \dot{H}^{s})\times
C([0,T); \dot{H}^{s-1})$, the necessary requirements to set up a fixed
point lead to a smaller Banach space. For example, this translates
into additional space-time integrability conditions, like $u\in
L^p_t(L^q_x)$ for suitable $p,q$. The resulting well-posedness result
is often deemed conditional (to these additional requirements). Our aim is to remove these
assumptions which are incorporated in the uniqueness and existence
statement given by, say, Picard's theorem, and prove unconditional
well-posedness (sometimes called unconditional uniqueness in the
literature), that is uniqueness in the natural class, where the flow
is continuous. Note that in the wave map
situation, one does not construct a solution directly by iteration, at
least when working at the critical regularity. Nevertheless, in order
to obtain a priori estimates, one is led to add similar requirements
($\partial u\in L^2_t(L^{2d}_x)$ for example in \cite{SS,NSU}).

From now on, we generically denote $(\nabla u,\partial_t u)$ as $\partial u$,
so that any statement regarding $u$ and $\partial_t u$ can be
summarized into one, like $\partial u\in C_t(\dot{H}^{s-1})$. Also, for any Banach space $X$, 
 $C(X) $ will denote the space $ C([0,T); X) $, $L^p(X)$ will denote
 the space $L^p(0,T;X)$ and $L^p(L^q)$  will denote the space 
$L^p(0,T; L^q (\mathbb{R}^d  ))$.
\begin{theorem}\label{main}
  Let $u$ be a solution to (\ref{wm}) on $[0, T^*)$,
   with $d\geq 4$. Then $u$ is the unique solution of  (\ref{wm}) in the class
$$
\partial u\in C(\dot{H}^{\frac{d}{2}-1}).
$$
\end{theorem}
\begin{remark}
The same result should hold for $d=3$ if one is willing to consider $
\partial u\in L^\infty(\dot{H}^{\frac{1}{2}+\e}),$ $\e>0$. In fact, both
schemes of proof from \cite{SS,NSU} work in that framework, modulo
technicalities related to the low regularity (one has to take into
account the null form structure in the elliptic equation).
\end{remark}
 
Recently, there  were many works  proving unconditional well-posedness
for several hyperbolic systems (see \cite{fabuni,MP06} for 
the critical wave equation, \cite{FT03,FPT03} for the nonlinear {S}chr\"odinger equation, 
\cite{MN09fe} for the {Z}akharov system, 
\cite{MN03cmp} for the {M}axwell-{D}irac system).
 There is now a huge literature about  unconditional well-posedness 
for parabolic equations such as the Navier-Stokes system
\cite{FLT00,LM01cpde,Masmoudi03nantes}.  

A desirable goal would be to prove that strong solutions to the wave
map equation do coincide with weak solutions in (spatial) dimension
$d=2$, at least before possible blow-up; recall that when $d=2$ the
scale-invariant space is the natural energy space. Such a goal
appeared totally out of reach when the present work was
started. However, a great deal of progress was made in recent years on
the strong Cauchy theory for $d=2$, eventually leading to global
well-posedness for large finite energy data
(\cite{ST1,ST2} for negatively curved compact targets, \cite{KS} for
the hyperbolic space $\mathbb{H}^2$ as target, and \cite{Tao} and references therein
again for the hyperbolic space). We hope our present high dimensional result
provides a clear view of uniqueness issues in a relatively
straightforward functional setting, but that the strategy itself will
be of interest in the more intricate lower dimensional setting.

\section{Existence of solutions : From large data to small data}
We recall  that results on global well-posedness at the
critical level for small data can actually be extended to local
well-posedness for large data. All these results rely on a gauge
change, which requires a certain connection (associated to the map) to
be small. In high dimension, one may choose the Coulomb gauge, and in the large data case, there is no reason for this
connection to be small, as the elliptic system linking the connection
to the map needs not  have a unique solution. Fortunately, one can
take advantage of a fundamental property of the wave equation,
namely the  finite
speed of propagation. The equation (\ref{wm}) can be seen as a
semilinear wave equation, and the nonlinearity is a local one, since  it
can be written  as a product of $\partial u$ and (function of) $u$. Remark that
after performing the gauge transform alluded to above, this local character is
lost, the new nonlinearity involves pseudodifferential operators.

Given an arbitrary initial data $(u_0, u_1)$, let us explain how we can
construct a local solution $u$ using the known results for small data.
We can choose $r>0$ small enough such that
$$  Sup_{x \in \R^d}  || (\partial u_0, u_1) ||_{\dot H^{{d \over 2 } -1 }(B(x,r))}
\leq  \eps_0   $$
 where $\eps_0$ is a small parameter which will be chosen later.
For each ball $B(x,r)$, the initial data $(u_0, u_1)$ can be extended
to the whole space by $(\tilde u_0, \tilde  u_1)$ in such a way that
$$   || (\partial  \tilde u_0,  \tilde  u_1) ||_{\dot H^{{d \over 2 } -1 }( \mathbb{R}^{d+1})}
\leq  C  \eps_0   $$
for some constant $C$ which only depends on $N$ and $r$.
For each $x \in  \mathbb{R}^d$, we  can use the results of
  \cite{SS} to  construct a global solution $u_x $ which is unique
  in the class  $\partial u  \in C(\dot H^{{d \over 2 } -1 } ( \mathbb{R}^{d+} ) ) \cap
L^2(L^{2d}_x) $. For each $t < { r \over 2}$, we can define $u$ by
$$  u(t,y) = u_x(t,y  )   \quad \hbox{if}  \quad (t,y)  \in  C(r,x) $$
where $C(r,x)$ denotes the backward light cone of vertex $(r,x)$
$$ C(r,x) = \{  (t,y)  \ |   \,    |y-x| \leq r - t \ \} .$$
We have only to make sure that  if $|x - x'| < 2r$ then
for all $(t,y)  \in  C(r,x)  \cap  C(r,x') $,  $u_{x'}(t,y) = u_x(t,y  )$.
Let $x_m = { x + x' \over 2 } $ and $r_m = r - |x-x_m|$ then
$C (r_m, x_m) =    C(r,x)  \cap  C(r,x') $.
Writing the energy estimate on  $w = u_x - u_{x'}$ in the cone
$C(r_m, x_m ) $ and using the same computation as in the
 uniqueness result of \cite{SS} (which uses the smallness condition
as well as the extra bound $L^2(L^{2d}_x) $), we infer that
$u_x =  u_{x'}$ in $C(r_m, x_m) $.

\section{Proof of theorem  \ref{main}}
 In order to avoid distracting dependencies on the dimension, we
shall restrict to the case $d=4$. The proof proceeds through several reductions.
We start with two solutions   $u$ and $v$   of (\ref{wm})
with the same initial data $(u_0, u_1)$ defined on some
time interval $[0, T)$ and such that
$\partial u, \,  \partial v \, \in C([0,T) ; \dot H^1)$.
 Without loss of generality, we can assume
that $u$ is the solution which was obtained in the previous
section (see also Shatah and Struwe
\cite{SS}). This solution $u$ is known to satisfy some
extra estimates which will be useful in the proof.

Next, notice that to prove uniqueness,
it is enough to prove that $u$ and $v$ coincide
  on some small time interval.
Indeed, if we can prove that   there exists $\tau$,    $0< \tau < T$ such that,
$\forall t, \, 0< t< \tau$, we have $u(t) = v(t)$ 
then by continuity, we deduce that $u(\tau) = v(\tau)$
and we can iterate the argument to prove that
$\forall t, \, 0< t< T$, we have $u(t) = v(t)$.

Finally, to prove uniqueness it is sufficient to
prove that $u$ and $v$ are equal on each backward light cone
with a vertex $(t,x)$ such that $0< t < r$.

\subsection{Gauge transform}
Given a (now small) data $\partial u_0$ and the wave map
$$
D_\alpha \partial^\alpha u=0,
$$
one can choose to work within the Coulomb gauge and take advantage of
carefully chosen frames to obtain a new system of the form
\begin{align}
  \label{eq:red}
  \Box q= & A \cdot \partial q+ q \partial\cdot A+A^2 q+ q (R(u)q^2) \\
 \label{eq:red2} \Delta A =  & \nabla (A^2)+\nabla (R(u)q^2),
\end{align}
where we simplified the  system  to a model case where $q$ is scalar,
$A$ is a vector and powers of $A$ are to be understood as bilinear
forms of its coefficients. To get a sense of perspective, one should
see $q\approx\partial u$, or more accurately any of the components of
the 1-form $du$, and $R$ should be seen as the curvature of the target
manifold, which we assume bounded along with all its derivative
(target with ``bounded geometry'').

For targets which are symmetric spaces, the coefficient $R(u)$ just
disappears. This new system can then be solved using an iteration
scheme, using Strichartz estimates up to the end-point (thus, the
restriction on $n\geq 4$), as is done in \cite{NSU}.

If one takes two small data which are the same, $\partial u_0=\partial
v_0$, the reduced system for $u$ and $v$ will be the same, and in
particular their respective data coincide. Hence, all there is to do
is to actually prove uniqueness for the system (\ref{eq:red}),(\ref{eq:red2}).

\subsection{A model case}
\label{model}
As we just saw, if the target happens to be a symmetric space, the renormalized wave map
system reduces to the following simpler system:
\begin{equation}
  \label{eq:easy}
\tag{RWM}
\left \{ \begin{array}{rcl}
  \Box q= & A \cdot \nabla q+ q \nabla\cdot A+A^2 q+q^3 \\
\Delta A =  & \nabla (A^2)+\nabla (q^2).
\end{array} \right.
\end{equation}
Essentially, the curvature term has disappeared in the elliptic
equation (which we refer to as (RWMe) while the wave equation part
will be (RWMh)), and we are left with a system involving only $q$ and $A$. Moreover,  if
we  make a smallness assumption, $A$ is entirely determined by
$q$, and in the present situation where (\ref{eq:easy}) is derived
from (\ref{wm}), we are under such an assumption.
\begin{theorem}
  The system  (\ref{eq:easy})  has a unique small solution in the class
  $\partial q\in L^\infty(L^2)$.
\end{theorem}
{\sl Proof.} Recall that one can perform a fixed point
in the class $E=C(\dot{H}^1)\cap L^2(\dot{B}^{1/6,2}_6)$ for $q$,
and $F=C(\dot{H}^1)\cap L^1(\dot{B}^{1,1}_4)$ for $A$
(\cite{NSU}). We therefore can prove uniqueness by comparing any
solution $v$ such that $\partial v\in C(L^2)$ with the reference
solution $u\in E$ (and its associated $A\in F$).
\begin{remark}
  Note that $\tilde F=L^\infty(L^4)\cap L^1(L^\infty)$ is enough to
  solve (and this is what happens, mutatis mutandis, in \cite{SS}),
  but the additional regularity information we carry in $F$ will be
  useful later.
\end{remark}
 The idea to
obtain uniqueness is to write an estimate at a lower regularity
level. This idea is recurrent when proving uniqueness for hyperbolic
systems, since taking differences yields a loss of derivative. In
fact, in \cite{SS}, uniqueness for $u\in C(\dot H^2)\cap L^2(L^8)$ is
established in this way, writing a difference estimate in $\dot
H^1$. At the level of $q$, this translate to uniqueness for $q\in
C(\dot H^1)\cap L^2(\dot W^{-1}_8)$ (though writing directly the
estimate at the $q$ level is most likely more involved than directly
on the true system).

 Consider $\delta=q-q'$ the difference between two solutions,
and $\alpha=A-A'$ the difference between the vectors, and
set $(q,A)$ to be the fixed point solution, namely the solution 
in $E\times F$. The equation for $(\delta,\alpha)$ will
be
\begin{equation}
  \label{eq:del}
  \tag{$\Delta$ RWM}
\left \{ \begin{array}{rcl}
  \Box \delta& \equiv & A \cdot \nabla \delta+\alpha
  \nabla(q-\delta)+\delta \nabla A+(q-\delta)\nabla
  \alpha\\
  & & {}+q\alpha(2A-\alpha)+\delta(A-\alpha)^2+\delta(q^2+q\delta+\delta^2) \\
\Delta \alpha &\equiv  & \nabla (2A\alpha-\alpha^2)+\nabla(2q\delta-\delta^2).
\end{array} \right.
\end{equation}

\begin{remark}
  We are in a situation where $q\in L^\infty(\dot H^1)$, small, say
  $\lesssim \varepsilon_0$. From
$$
A=|\nabla|^{-1}(A^2+q^2),
$$
we know that we can solve this elliptic equation:
$$
\|A\|_{\dot{H}^1}\lesssim \|A\|^2_{L^4}+\varepsilon_0^2,
$$
hence
$$
\|A\|_{L^4}\lesssim \|A\|^2_{L^4}+\varepsilon_0^2,
$$
which gives both $\|A\|_{L^4}$ and $\|A\|_{\dot{H}^1}$ small, and the
same is true also for $q'$, $A'$. We also
have space-time estimates on $A$ from what we know on $q$: say $q\in
L^2(\dot B^{1/6,2}_6)$, then $q^2\in L^1(\dot B^{0,1}_4)$, and we can
write
$$
\|A\|_{L^1(\dot B^{1,1}_4)}\lesssim \|A\|_{L^\infty \dot H^1}\|A\|_{L^1(\dot B^{1,1}_4)}+\|q^2\|_{ L^1(\dot B^{0,1}_4)}.
$$
Such an estimate will prove useful later, note that this immediately
gives $A\in L^1L^\infty$.
\end{remark}
We will write an estimate for $\delta$ in the following Strichartz space 
$$
X=C(\dot{H}^{\frac 1 6})\cap L^2(\dot{B}^{-\frac{2}{3},2}_6).
$$
For $\alpha$, one may think that, heuristically,
$\alpha=|\nabla|^{-1}(\delta^2)$, and this leads to (using one factor
$\delta$ in $C(\dot H^1)$ and the other factor $\delta$ in $
L^2(\dot{B}^{-\frac{2}{3},2}_6)$)
$$
\alpha\in Z=L^2(\dot B ^{1,2}_{12/7})\hookrightarrow L^2(L^3).
$$
\begin{remark}
  Let us motivate the choice of $X$ and $Z$: if it was not for the
  derivatives, the model equation for $q$ would be $\Box q=q^3$. In
  \cite{fabuni}, uniqueness for this equation is established for $\dot
  H^1$ data; this relies on a contraction estimate, with
  $\delta\in X_{-\frac 1 3}$ where $X_s=L^2(\dot{B}^{-s,2}_6)$. In fact, one has
  some freedom in the choice of $X$, and any $X_{s}$ with $-1<s\leq
  -\frac 1 3$ would do. However, from the embedding $\dot
  H^1\hookrightarrow \dot{B}^{-\frac{1}{3},2}_6$ the choice $X_{-\frac
  1 3}$ seems straightforward. In our setting, however, the source
  term is more like $|\nabla|^{-1}(\delta^2) \nabla \delta$. Using our
  knowledge $\nabla \delta \in L^2$, the requirement on $|\nabla|^{-1}(\delta^2)$ becomes
  $L^2(\dot{B}^{s+\frac 2 3}_6)$; Part of the product in the source involves low frequencies of this term, producing the worst possible situation. Since we cannot afford to use Sobolev
  embedding, this requires $s+\frac 2 3\leq 0$. Thus one
  is led naturally to pick $s=-\frac 2 3$. This in turn requires to
  check that $\delta$ belongs to the chosen $X$ space, which we do below.
\end{remark}
\begin{lemma}
  If $u,v$ are two solutions of (\ref{eq:easy}) such that $\partial
  u,\partial v\in C(L^2)$, then $\delta \in X$ and $\alpha\in Z$.
\end{lemma}
Proof. Certainly $\delta,\alpha,q,A \in \dot{H}^1\hookrightarrow
L^4 \hookrightarrow \dot{B}^{-1/3,2}_6$. Using the equation and
looking only at the worst possible term,
$$
|\nabla|^{-1}(ab)\nabla c \in L^{\frac 4 3}\hookrightarrow \dot{H}^{-1},
$$
where $a,b,c$ are $\delta,\alpha,q,A$, and $ab\in (L^4)^2=L^2$, hence
$|\nabla|^{-1}(ab)\in \dot H^1\hookrightarrow L^4$ and $\nabla c\in
L^2$.

By Strichartz estimates, we deduce that  $\delta\in C(\dot{H}^{0})\cap
L^2(\dot{B}^{-\frac{5}{6},2}_6)$ and by interpolation (recall we are
local in time, $L^2_t$ is controlled by $L^\infty_t$) we obtain
$\delta\in X$.
\begin{remark}
  In all the remaining part of the paper, we will have to perform
  various product estimates, and rely heavily on the first part of the
  Appendix to do so. We refer to the Appendix for precise definitions
  of the LF/MF/HF interactions, in connection with paraproduct
  decomposition. We simply recall here that LF (resp. MF,HF) is
  meant for low frequencies (resp. medium, high) interactions between
  frequencies of factors in a product. 
\end{remark}
Now we check that similarly $\alpha\in Z$: using the elliptic equation
again, one has to check that $\delta q\in L^2 (\dot B^0_{12/7})$. For
this, we use Proposition \ref{prodpar1}: the
MF-HF interaction is easily treated, as $\delta \in L^2(\dot B^{-2/3}_6)$ and $q \in
L^\infty(\dot B ^1_2)$, therefore by \eqref{eq:MFHF} this interaction will be in
$L^2(\dot B^{1/3}_{3/2})\subset L^2(\dot B^0_{12/7})$; the HF-LF interaction is just as fine since 
$\delta \in L^\infty(\dot
B^{1/6}_2) $ and $q \in L^2(\dot B^{1/6}_6)\subset L^2(\dot B
^{-1/6}_{12})$ which results in a $L^2(\dot
B^0_{12/7})$ term from \eqref{eq:MFHF}. We have
$$
\|\alpha\|_Z\lesssim \varepsilon_0 \|\alpha\|_Z+\varepsilon_0
\|\delta\|_X+\|\delta\|_X\|q\|_{L^2 (\dot B^{1/6}_6)}\,,
$$
and the $\delta^2$ is just as easy, HF in $L^\infty (\dot H^1)$ and MF
in $L^2(B^{-2/3}_6)$. This ends the proof of the lemma.

 We now aim at closing an estimate in $X$. We will prove
 \begin{proposition}
   Let $\delta=u-v$ be the difference of two solutions with the same
   initial data. Then
   \begin{equation}
     \label{eq:contract}
     \|\delta\|_X\lesssim \varepsilon_0 \|\delta\|_X,
   \end{equation}
from which we can infer $\delta=0$.
 \end{proposition}
Proof. The source term in the
 equation on $\delta$ should be anywhere (in term of interpolation)
 between $L^2(L^{6/5})$ and $L^1(\dot H^{-5/6})$, which are then pulled 
 back to $X$ by $\Box^{-1}$. Indeed, let us   denote $X'=L^2(L^{6/5})+L^1(\dot H^{-5/6})$ and 
recall the following end-point Strichartz estimate: let $\Box \delta =
F$ and $  (\delta, \partial_t \delta ) (t=0) =0 $, then
\begin{equation}
  \label{eq:eps}
 \| \delta \|_X  \lesssim  \| F \|_{X'}.
\end{equation}
Let us go through each term  appearing in the right-hand side of ({$\Delta$ RWM}).
\begin{itemize}
\item The easy terms: $\alpha \nabla (q- \delta),$ $(\nabla \alpha)
  (q-\delta) $.
   These can all be
  dealt with by  Sobolev embedding and H\"older inequality.

 As $\alpha\in L^2(L^3)$ and
  $\nabla\delta,\, \nabla q\in L_t^\infty L^2)$ we have the corresponding
  term in $L^2(L^{6/5})$.
Similarly, $\alpha\in L^2(\dot B^{1,2}_{12/7})$ and $q-\delta \in
    L^\infty( L^4)$ give $L^2(L^{6/5})$.
Hence, we have
\begin{equation}\label{est1}
    || \alpha \nabla(q-\delta)  +  (q-\delta)  \nabla \alpha ||_{X'}
     \lesssim  || q - \delta ||_{E} ||\alpha ||_{Z}\lesssim
     \varepsilon_0 \|\delta\|_X.
\end{equation}
\item The term $\delta^3$: interpolation between $\delta\in C(\dot
  H^1)$ and $\delta\in L^2(\dot B^{-2/3}_6)$ yields
  $\delta\in L^4(\dot B^{1/6,2}_{3})\hookrightarrow L^4(L^{24/7})$, hence
  $\delta^2\in L^2(L^{12/7})$. Then using that  $\delta\in L^\infty (L^4)$, one gets
  that  $ \delta^3 \in L^2(L^{6/5})$, namely 
$$
\| \delta^3\|_{X'}\lesssim \|\delta\|^2_E \|\delta\|_X.
$$
\item The two terms $A\nabla \delta$ and $\delta \nabla A$: here, one
  has to perform a paraproduct decomposition and deal with the different
  frequency interactions in a suitable way. We follow the conventions set up in the
  Appendix for the different interactions in the paraproduct
  decomposition of a product.
  \begin{itemize}
  \item LF-HF interaction: $\nabla \delta \in C(\dot H^{-5/6})$, hence we are
  forced to have $A\in L^1(L^\infty)$, which is fortunately true since
  $A$ is the good connection from local Cauchy theory (recall actually $A\in L^1(\dot B^{1,1}_4)$).
\item HF-MF interaction: same information on $\delta$, but using $A\in L^1(\dot B^{1,1}_4)$ (so that
  $1+-5/6>0$) and embedding. Notice how we need the regularity on $A$
  ($A\in L^1(L^\infty)$ would be too weak).
  \end{itemize}
We thus have
$$
\|A\nabla \delta\|_{X'}\lesssim \|A\|_F \|\delta\|_X.
$$
For $(\nabla A) \delta$, we proceed similarly:
  \begin{itemize}
  \item LF-HF interaction: $\delta\in C(\dot H^{1/6})$ and $\nabla
    A\in L^1(\dot B^{-1,1}_\infty)$ yields a term in $L^1(\dot H^{-5/6})$.
\item HF-LF interaction: similarly, $\nabla A\in L^1(L^4) $,
  $\delta\in C(L^{24/11})$ yields a term in $L^1(\dot H^{-5/6})$ after embedding.
\item HF-HF interaction: again,  $\nabla A\in L^1(L^4)$ and $\delta\in C(\dot H^{1/6})$
  yields a term in  $L^1(\dot H^{-5/6})$ after embedding. 
  \end{itemize}
Thus, we obtain
$$
\|\delta \nabla A\|_{X'}\lesssim \|A\|_F \|\delta\|_X.
$$
\item The terms $\delta^2q$ and $\delta q^2$; we only proceed with the
  details of $\delta q^2$, $\delta^2 q$ being easier:
  $ q^2$ is just like $\nabla A$, that is $q^2\in L^1(L^4)$, hence we do
  exactly as the previous one $\delta \nabla A$.
  \end{itemize}
Summing all the above estimates, we get
$$
\|\delta\|_X\lesssim \|\text{R.H.S}\|_{X'} \lesssim \varepsilon_0 \|\delta\|_X,
$$
where $\text{R.H.S.}$ denotes the source term in (\ref{eq:del}), which
ends the proof.
\subsection{The general case}

We start with two solutions $u$ and $v$ such that $\partial u $
and $\partial v$ are in $C((\dot H^1)$. Without loss of
generality, we can assume that $u$ is the solution constructed
by Shatah and Struwe in \cite{SS}.  The argument given
in the last section was based on some smallness condition.
Using the finite speed of propagation for the wave
equation we will reduce our problem to the small case.
We choose $r > 0$ small enough that
$$  Sup_{x \in \R^d}  || \partial u(0) ||_{\dot H^1(B(x,r))}
\leq  \eps_0   $$
 where $\eps_0$ is a small parameter which will be chosen later.

Next, using the continuity of $u$ and $v$ with respect to time,
we can choose $\tau$,   $ 0 < \tau \leq {r \over 2 }   $ such that
$$  \forall  t, \  0 \leq t \leq \tau, \
     || \partial u (t) -  \partial  u(0 ) ||_{ \dot H^1}
  +  || \partial v (t) -  \partial  u(0 ) ||_{ \dot H^1}  \leq \eps_0.  $$

To prove that $u$ and $v$ coincide for all $0\leq t < \tau$, 
it is sufficient to prove that they coincide  on each truncated
backward light cone of the form
 $$ C^\tau_0(x_0,r) =  \{ (t,x) \  / \  |x - x_0 | \leq r - t \ , \
   0 \leq t \leq \tau \} .  $$
In the sequel, we restrict ourselves to the  backward light cone of
center $x_0 = 0 $.
The proof will be reduced to the proof given in the model
case. There are only two extra  difficulties we have to
handle.

The first difficulty lies in the choice of some coordinate
system where we can write our equation in a form similar to
$(RWM)$. As in \cite{SS}, we have to choose a frame $e$
and a connection $A$ satisfying the Coulomb gauge and
express $\partial u $ in that frame using the
coordinates $q$.  Moreover, to get good estimates
for $\alpha = A - A'$ and $e-e'$ in terms of
$\delta= q - q'$ we have to  construct    local  frames.
The second difficulty comes from the fact that the
curvature tensor is no longer constant
and an extra  term $R(u)$ will appear in the equation
$(RWM)$ and hence we have to estimate
$R(u) - R(v) $ in terms of $\delta$.

Let us start  by constructing the the local frame
$e$ and $e'$ associated  respectively  to $u$ and $v$.

\subsubsection{Construction of the frame}

Without loss of generality, we can assume that
$TN$ is parallelizable; hence, we can find  smooth  vector fields
$\bar e_1, ..., \bar e_k$ such that at each $p \in
N$ the family $\{ \bar e_1, ..., \bar e_k \}$ is an
orthonormal basis of $T_p N$ (see for instance \cite{Helein}). Given the
map $u$ or $v$ from $R^{d + 1 }$ into $N$, the
family  $\{ \bar e_1 \circ  u  , ..., \bar e_k \circ  u  \}$
is a smooth orthonormal frame of the pull-back bundle
$u^* TN$. Moreover, we may freely rotate this frame
at any point $(t,x) \in \R^{d+1} $ with a matrix
$(R_a^b) = (R_a^b(z)) \in  SO(k) $, thus obtaining the frame

\begin{eqnarray}
  \label{eq:change}
  e_a =  R^b_a \bar  e_b \circ u \,   , \quad \quad  1 \leq a \leq k.
\end{eqnarray}
For our uniqueness proof, it will be important that  $R^b_a  $ only depend 
on the solution in the cone $ C^\tau_0(x_0,r) $. 
We choose $R^b_a$ in the truncated cone $C_\tau$ by minimizing
for each time $0 \leq t < \tau$ the following functional

\begin{equation}
  \label{eq:minimize}
  F(R) =\int_{B(0, r - t )} \sum_{i=1}^{d} \sum_{a,b =1}^{k}
   \Big<   \frac{\partial e_a}{\partial x^i} , e_b\,.
  \Big>^2  \ dx
\end{equation}

The existence of a minimizer can be proved
following the same proof as in \cite{Helein}.
Moreover, denoting $A^a_{b,\alpha}
 = \Big<   \frac{\partial e_a}{\partial x^\alpha} , e_b
  \Big> $ for $ 1 \leq \alpha \leq d $ and $1\leq a,b \leq k$,
we get the following Euler-Lagrange equation
\begin{equation}
  \label{eq:euler-lagrange}
\left\{ \begin{array}{cccc}
  \partial_\alpha A^a_{b,\alpha}  &= 0  \quad & \hbox{in} & \quad
      B(0, r - t )         \\
      A^a_{b,\alpha}. n_{\alpha} &= 0   \quad & \hbox{on} & \quad
     \partial B(0, r - t )
\end{array} \right.
\end{equation}
where $n$ denotes the normal to the ball $B(0, a - t )  $.
We need an extra equation to determine  $A$, which we can get
from the curvature of the pull-back covariant derivative
$D = (D_\alpha)_{0\leq \alpha \leq d}  $.
Indeed, using that the Lie bracket between  $D_\alpha$
and $D_\beta$ vanishes, we get that
for all  $ 1 \leq  a , b \leq k  $
\begin{equation}
\label{eq:cuvature}
  \partial_\alpha  A^a_{b,\beta} -
  \partial_\beta  A^a_{b,\alpha}
  + [ A_{\alpha} , A_{\beta} ]^a_b  = R(u) ( \partial_\alpha u,
  \partial_\beta  u )\,.
\end{equation}
Expressing $du$ in the frame $e$, we get $\partial_\alpha u =
q^a_\alpha e_a$. On the other hand, written in the $q$ coordinate
system, the wave map equation yields
\begin{equation} 
  \label{eq:wave-map-q}
   \Box q_\beta     \equiv  2 A^\alpha
 \partial_\alpha q_\beta  +   ( \partial^\alpha  A_\alpha )
q_\beta   + A^\alpha A_\alpha q_\beta + F^\alpha_\beta q_\alpha.
\end{equation}
Notice that the system of equations we obtained,
namely (\ref{eq:wave-map-q}), (\ref{eq:euler-lagrange}) and
  (\ref{eq:cuvature})   is very similar to the model
problem we studied in the previous subsection.

The construction of the frame $e$, the connection $A$ and
the components $q$ of $du$ can be carried out also for
the solution $v$. We denote  $e'$, $A'$ and $q'$ the
resulting frame, connection and components of $dv$.
In the sequel, we denote $\delta = q-q'$,
$\alpha = A - A'$.
Using the different equations we have at hand,
we will estimate  $du-dv$ and $\alpha$ in terms
of $\delta$ and then prove a closed estimate for
$\delta$ from which we deduce that $\delta $ should vanish
as well as $ du-dv$ and $\alpha$.
\subsubsection{Proof of Theorem \ref{main}}
In the generic situation, one has (as a model) the following two
equations: the first one is the wave equation \eqref{eq:wave-map-q} holding inside a
space-time cone, and the second is an    elliptic equation. It  is a short-hand
for what is really an elliptic div-curl system (\eqref{eq:euler-lagrange} and
\eqref{eq:cuvature}) holding on fixed time slices with appropriate
boundary conditions. The elliptic theory yields the same regularity
estimates in that case as in our simplified model with a Laplacian.
\begin{equation}
  \label{eq:generic}
\tag{RWM}
\left \{ \begin{array}{rcl}
  \Box q\equiv & A \cdot \nabla q+ q \nabla\cdot A+q^3 \\
\Delta A\equiv  & \nabla (A^2)+\nabla (R(u)(q^2)).
\end{array} \right.
\end{equation}
$R$ is essentially a smooth function related to the curvature tensor
on the target space (hence connected with the Christoffel symbols associated with the
connection), which we assume bounded and with all derivatives bounded. Problems  may arise whenever we encounter $R(u)-R(u')$, however,
$$
R(u)-R(u')=(u-u')\int_0^1 R'(\theta u+(1-\theta)u')d\theta,
$$
Provided we seek an estimate on $R(u)-R(u')$ such that we only use
$R'\in L^\infty_{t,x}$,  we are left with the other factor, namely$(u-u')$. However, we can control $\partial(u-u')$ by $\delta$, and it
turns out to be sufficient to close the estimates. The new system for
$(\delta,\alpha)$ is (denoting by $w=u-u'$)
\begin{align}
  \label{eq:deldur}
  \Box \delta  \equiv & A \cdot \nabla \delta+\alpha
  \nabla(q-\delta)+\delta \nabla A+(q-\delta)\nabla
  \alpha \nonumber  \\
   & {}+\delta(q^2+q\delta+\delta^2)  \tag{$\delta$ RWM}\\
\alpha  \equiv  & |\nabla|^{-1}( (2A\alpha-\alpha^2)+R(u-\delta)(2q\delta-\delta^2)+(R(u)-R(u-w))q^2).\nonumber
\end{align}
Thus the modification appears in the elliptic equation on the
connection. From the computation in the previous section, one infers
that $2q\delta-\delta^2\in L^2(L^{\frac {12} 7})$, and combined with
$R\in L^\infty(L^\infty)$ we dispose of the first term with $R$ as we did
in the model case. The
next term is the real novelty here. Assuming that $\partial w \equiv
\delta$, we get that
$$
w\in L^2(\dot B^{\frac 1 3}_6) \hookrightarrow L^2(L^{12}),
$$
and using $q^2\in L^\infty (L^2)$, we get the desired $L^2 (L^{\frac {12} 7})$ estimate for the source term.

All other terms may be estimated like in the model case, and we can
therefore close an estimate on $\delta$ as we did in the previous
section, up to localization to the interior of the light cone for
space-time estimates and localization to balls for elliptic
estimates. Fortunately, all the required estimates may easily be
transposed in such a situation, as explained in Appendix
\ref{sec:albs}. This ends the proof.

\section*{Acknowledgments}

The first author was partially supported by
 NFS grant DMS-0703145. The second author wishes to thank the Courant
 institute of mathematical science for its kind hospitality during a
 visit where this work was initiated; he was partially supported by A.N.R. grant SWAP.

\appendix
\appendixpage

\section{Product estimates}
\label{sec:ape}
In this appendix we describe various product estimates in Besov
spaces, which are used throughout the rest of the paper. We do not
claim novelty here, but we do however emphasize that thinking about
the product in terms of different frequency interactions is crucial
in our situation.
\begin{proposition}
\label{prodpar1}
  Let $f\in \dot{B}^{s_1,q_1}_{p_1}=B_1$ and $g\in
  \dot{B}^{s_2,q_2}_{p_2}=B_2$. Assume $s_i-\frac{d}{p_i}<0$, define $r_i$
  such that $s_i-\frac{d}{p_i}=-\frac{d}{r_i}$ (Sobolev embedding
  exponent if $s_i>0$), and assume moreover that
  $\frac{1}{r_1}+\frac{1}{r_2}<1$.
  \begin{enumerate}
  \item Suppose $s_1>0$ and $s_2<0$, and $r_1\geq q_1$. Then, $fg=\pi_1+\pi_2$ where $\pi_1\in
    \dot{B}^{s_1+s_2,q}_{p}$, $\pi_2\in \dot{B}^{s_2,q_2}_{P_2}$, with
$$
\frac{1}{p}=\frac{1}{p_1}+\frac{1}{p_2}\,\,,\,\frac{1}{q}=\frac{1}{q_1}+\frac{1}{q_2}\,\,,\,
\frac{1}{P_2}=\frac{1}{p_2}+\frac{1}{r_1},
$$
and
\begin{equation}
  \label{eq:MFHF}
  \|\pi_1\|_{\dot{B}^{s_1+s_2,q}_{p}}+\|\pi_2\|_{\dot{B}^{s_2,q_2}_{P_2}}
  \lesssim \|f \|_{B_1}\|g\|_{B_2}.
\end{equation}
We  call $\pi_1$ the high-medium frequencies interaction term, the high
frequencies referring to the $f$ factor and the medium to the $g$ factor. We abbreviate it to
$HF-MF$ (or $MF-HF$ if $f$ and $g$ are switched). Similarly, $\pi_2$ is the
low-high frequencies interaction term, or $LF-HF$ for short.
\item Suppose $s_1,s_2>0$ and $r_i\geq q_i$, then  $fg=\pi_1+\pi_2+\pi_3$ where $\pi_3\in
    \dot{B}^{s_1+s_2,q}_{p}$, $\pi_1\in \dot{B}^{s_1,q_1}_{P_1}$ $\pi_2\in
    \dot{B}^{s_2,q_2}_{P_2}$, with $p,q,P_2$ as above,
    $\frac{1}{P_1}=\frac{1}{p_1}+\frac{1}{r_2}$ and
\begin{equation}
  \label{eq:HFHF}
  \|\pi_3\|_{\dot{B}^{s_1+s_2,q}_{p}}+\|\pi_1\|_{\dot{B}^{s_1,q_1}_{P_1}}+\|\pi_2\|_{\dot{B}^{s_2,q_2}_{P_2}}
  \lesssim \|f \|_{B_1}\|g\|_{B_2}.
\end{equation}
We refer to $\pi_1$ as the $HF-LF$ term, $\pi_2$ as the $LF-HF$ term
and $\pi_3$ as the $HF-HF$ term, similarly to the previous case.
  \end{enumerate}
\end{proposition}
Such product estimates are classical, see
e.g. \cite{RS}. Consider the first case: we decompose $fg$ as
$$
fg=\pi_1+\pi_2=\sum_j S_{j+2} g \Delta_j f +\sum_j  S_{j-1} f \Delta_j g.
$$
The term $\pi_2$ is a sum of frequency localized pieces, meaning that
for a finite number of $k$ close to $j$,
$$
\Delta_j \pi_2=\sum_{k\approx j} \Delta_j( S_{k-1} f \Delta_k g).
$$
For convenience we only deal with the $k=j$ term. For the low
frequencies $S_{j-2} f$, we use Sobolev embedding  and the fact that $r_1\geq q_1$ to get
$$
\|S_{j-2} f\|_{r_1}\lesssim \|f\|_{B_1}.
$$
For the high frequencies $\Delta_j g$,
$$
\|\Delta_j g\|_{p_2}\lesssim 2^{s_2 j} \varepsilon_j \|g\|_{B_2},
$$
where $\varepsilon_j \in l^{q_2}$. The result follows by H\"older.

The other term $\pi_1$ is a sum of dyadic terms localized in balls of
radius $2^j$. We estimate
$$
\Delta_j \pi_1 =\sum_{j\lesssim k} \Delta_j(\Delta_k f S_{k+2} g),
$$
and, since $s_2<0$ and recalling $S_j=\sum_{l<j}\Delta_l$,
$$
\|S_{j+2} g\|_{p_2}\lesssim 2^{-s_2 j} \mu_j \|g\|_{B_2},
$$
with $\mu_j \in l^{q_2}$. Thus
$$
\|\Delta_j \pi_1\|_{p}\lesssim \sum_{j\lesssim k}  2^{-(s_1+s_2) k}
\mu_k \eta_k \|g\|_{B_2}\|f\|_{B_1}= 2^{-(s_1+s_2) j}\lambda_j
\|g\|_{B_2}\|f\|_{B_1},
$$
with $\lambda_j \in l^q$, and we are done.

The other case proceeds similarly, except we use the full paraproduct
decomposition, namely
$$
fg=\pi_1+\pi_2+\pi_3=\sum_j  S_{j-1} g \Delta_j f+\sum_j  S_{j-1} f \Delta_j
g+\sum_{|k-k'|\leq 2}\Delta_k f\Delta_{k'}g.
$$
The first two terms are treated like the term $\pi_2$,
and the term $\pi_3$ is treated as the term $\pi_1$.

\section{Localizing Besov spaces}
\label{sec:albs}
In our context and in order to take advantage of the finite speed of
propagation, one wishes to localize all the usual estimates to a
backward light cone. Such a procedure is well-known in the context of
the critical semilinear wave equation: see for example \cite{SS93},
where an explicit extension procedure is given to achieve this goal.

In our setting we do not need to worry about space-time Besov spaces
and the geometry of cones: we are always interested in estimates in a
truncated backward cone (avoiding the tip of the cone). Denote by
$(x_0,t_0)$ the vertex of such a backward cone $C(x_0,t_0)$ ($t_0>0$),
we consider a slab $C_0^a(x_0,t_0)=C(x_0,t_0)\cap \R^d\times
[0,a]$. For each $t\in [0,a]$, denote by $B_t$ the corresponding time
slice of $C_0^a(x_0,t_0)$: this (space) ball or radius $t$ is a smooth
domain of $\R^d$.
\begin{definition}[\cite{Tri}]
A function $f(x)\in \mathcal{D}'(B_t)$ belongs to the (spatial) Besov
space $\dot B ^{s,q}_p(B_t)$ iff there exists $g\in  \dot B
^{s,q}_p(\R^d)$ such that $f$ is the restriction of $g$ to $B_t$ (as
distributions). The norm of $f$ is then the minimum over all possible
extensions $g$ of their Besov norm in $\mathbb{R}^d$.
\end{definition}
The main property we need is the existence of an extension operator:
if we call $R$ the restriction operator, i.e. $f=Rg$, then there
exists an operator $E$ from $\dot B ^{s,q}_p(B_t)$ to $ \dot B
^{s,q}_p(\R^d)$ such that $f=REf$. Moreover, this extension operator
can be chosen to be the same whenever $(p,q,s)$ are in a bounded
domain (which is always the case for us). We refer again to \cite{Tri}
for a detailed presentation.

Next, for a given space-time function $u(x,t)$, we can define
$u\in L^p_T(\dot B ^{s,q}_p (B_t))$ by
$$
\int_0^T \|u(\cdot,t)\|^p_{\dot B ^{s,q}_p (B_t)}\,dt <+\infty.
$$

After these preliminaries, we can localize estimates in this way:
\begin{itemize}
\item Product estimates.\\
Let $f_1\in L^{r_1}_T \dot B^{s_1,q_1}_{p_1}(B_t),f_2 \in L^{r_2}_T
\dot B^{s_2,q_2}_{p_2}(B_t)$. Then $f_1 f_2\in  L^{r}_T \dot
B^{s,q}_{p}(B_t)$ and
$$
 \|f_1 f_2\|_{
  L^{r}_T \dot B^{s,q}_{p}(B_t)}\lesssim  \|f_1 \|_{ L^{r_1}_T \dot
  B^{s_1,q_1}_{p_1}(B_t)} \|f_2\|_{L^{r_2}_T
\dot B^{s_2,q_2}_{p_2}(B_t)},
$$
where $r,s,p,q$ are the same as in $\R^d$. In fact, we have $g_1$ and
$g_2$ the extensions of $f_1$ and $f_2$, we perform the product $g_1
g_2$ and then we have $f_1 f_2=R(g_1 g_2)$, and the inequality between
the two norms.
\item Linear estimates for the wave equation.\\
Consider first the inhomogeneous wave equation, $\Box u=F$, with zero
Cauchy data at time $t=0$. By finite speed of propagation, $u$ inside
the backward cone $C(t_0,x_0)$ depends only on $F$ inside the same
region. Moreover (causality), $u$ at $t=a$ depends only on $F$ on
$C_0^a(t_0,x_0)$. Given $F\in L^{r'}_a(\dot B^{s',2}_{p'}(B_t))$
where $(r',p')$ is a dual admissible Strichartz pair, we can extend it
to an $\tilde F\in   L^{r'}_a(\dot B^{s',2}_{p'}(\R^d))$, apply
Strichartz estimates, recover $\tilde u=\Box^{-1} \tilde F$ such that
$\tilde u\in   L^{\lambda}_a(\dot B^{s,2}_{\mu}(\R^d))$ where
$(\lambda,\mu)$ is any admissible Strichartz pair and $u=R\tilde
u$. Therefore,
$$
\|u\|_{ L^{\lambda}_a(\dot B^{s,2}_{\mu} (B_t))}\lesssim 
\|F\|_{L^{r'}_a(\dot B^{s',2}_{p'}(B_t))}.
$$ 
One can proceed similarly for the data to obtain the full range of estimates
for the Cauchy problem.
\end{itemize}
Combining these two observations, we can localize all the
estimates from Subsection \ref{model} without modification, whenever
we are facing a product of functions or an estimate on a
solution to the wave equation (through the use of the Duhamel
formula).

 \end{document}